\documentclass{gtart_h}

\def\ifplaintex{\expandafter\ifx\csname documentclass\endcsname\relax}

\def\gtp{{\mathsurround=0pt\it $\cal G\mskip-2mu$eometry \&\ 
$\cal T\!\!$opology $\cal P\!$ublications}}  

\def\recd{{\small Received:\qua\receiveddate\ifx\reviseddate\relax
\else\qquad Revised:\qua\reviseddate\fi\par}} 

\def\lognumber#1{\def\thelognumber{#1}}
\def\volumenumber#1{\def\thevolumenumber{#1}}
\def\volumeyear#1{\def\thevolumeyear{#1}}
\def\papernumber#1{\def\thepapernumber{#1}}
\def\pagenumbers#1#2{\def\startpage{#1}\def\finishpage{#2}}
\def\published#1{\def\publishdate{#1}}

\def\received#1{\def\receiveddate{#1}}
\def\revised#1{\def\reviseddate{#1}}
\def\accepted#1{\def\accepteddate{#1}}

\def\asciiaddress#1{\def\theasciiaddress{#1}}
\def\asciiemail#1{\def\theasciiemail{#1}}

\long\def\asciiabstract#1{\long\def\theasciiabstract{#1}}

\let\\\par\let\thelognumber\relax\let\thevolumenumber\relax
\let\thepapernumber\relax\let\thevolumeyear\relax\let\startpage\relax
\let\finishpage\relax\let\publishdate\relax\let\receiveddate\relax
\let\reviseddate\relax\let\accepteddate\relax\let\theasciititle\relax
\let\theasciiauthors\relax\let\theasciiaddress\relax
\let\theasciiabstract\relax

\let\theasciiemail\relax

\ifplaintex
\font\logobig=cmssbx10 scaled 3836
\font\logomed=cmssbx10 scaled 2557
\else
\font\logobig=cmssbx10 scaled 4200
\font\logomed=cmssbx10 scaled 2800
\fi

\long\def\makeagttitle{   
\count0=\startpage
\agt\hfill      
\hbox to 45truept{\vbox to 0pt{\vglue -13truept{\logomed A\kern -.37em{\logobig 
T}\kern -.38em G}\vss}\hss}
\break
{\small Volume \thevolumenumber\ (\thevolumeyear)
\startpage--\finishpage\nl
Published: \publishdate}

\vglue .25truein

{\parskip=0pt\leftskip 0pt plus
1fil\def\\{\par\smallskip}{\Large\bf\thetitle}\par\medskip} \vglue
0.05truein

{\parskip=0pt\leftskip 0pt plus 1fil\def\\{\par}{\sc\theauthors}
\par\medskip}%
 
\vglue 0.03truein 

{\small\leftskip 25truept\rightskip 25truept{\bf Abstract}\stdspace\theabstract

{\bf AMS Classification}\stdspace\theprimaryclass
\ifx\thesecondaryclass\relax\else; \thesecondaryclass\fi\par
{\bf Keywords}\stdspace \thekeywords\par}\vglue 7truept

}   

\ifplaintex
\font\phead=cmsl9 scaled 950
\font\pnum=cmbx10 scaled 913
\font\pfoot=cmsl9 scaled 950
\headline{\vbox to 0pt{\vskip -4.5mm\line{\small\phead\ifnum
\count0=\startpage ISSN 1472-2739 (on-line) 1472-2747 (printed)
\hfill {\pnum\folio}\else\ifodd\count0\def\\{ }%
\ifx\theshorttitle\relax\thetitle\else\theshorttitle\fi\hfill{\pnum\folio}
\else\def\\{ and }{\pnum\folio}\hfill\ifx\theshortauthors\relax\theauthors
\else\theshortauthors\fi\fi\fi}\vss}}
\footline{\vbox to 0pt{\vglue 0mm\line{\small\pfoot\ifnum\count0=\startpage
\copyright\ \gtp\hfill\else
\agt, Volume \thevolumenumber\ (\thevolumeyear)\hfill\fi}\vss}}
\else
\font=cmsl9 scaled 1050
\font\lnum=cmbx10 
\font=cmsl9 scaled 1050
\makeatletter
\def\@oddhead{{\small\lhead\ifnum\count0=\startpage ISSN 1472-2739 
(on-line) 1472-2747 (printed)\hfill {\lnum\number\count0}\else\ifodd\count0
\def\\{ }\ifx\theshorttitle\relax \thetitle \else\theshorttitle\fi\hfill
{\lnum\number\count0}\else\def\\{ and }{\lnum\number\count0}
\hfill\ifx\theshortauthors\relax 
\theauthors\else\theshortauthors\fi\fi\fi}}\def\@evenhead{\@oddhead}
\def\@oddfoot{\small\lfoot\ifnum\count0=\startpage\copyright\ \gtp\hfill\else
\agt, Volume \thevolumenumber\ (\thevolumeyear)\hfill\fi}
\def\@evenfoot{\@oddfoot}
\makeatother
\fi
\let\maketitlepage\makeagttitle
\let\makeshorttitle\maketitlepage
\let\maketitle\maketitlepage

\newwrite\gtoutfile
\long\gdef\makeheadfile{  
{\def\\{, }\def\s{ }
\immediate\openout\gtoutfile head.xxx
\immediate\write\gtoutfile{Proxy-for: \ifx\theasciiauthors\relax
\theauthors\else\theasciiauthors\fi\s<\ifx\theasciiemail\relax\theemail\else\theasciiemail\fi>}
\immediate\write\gtoutfile{\noexpand\\}
\immediate\write\gtoutfile{Authors: \ifx\theasciiauthors\relax
\theauthors\else\theasciiauthors\fi}
{\def\\{ }\immediate\write\gtoutfile{Title: \ifx\theasciititle\relax
\thetitle\else\theasciititle\fi}}
\immediate\write\gtoutfile{Subj-class: GT or SG, GR etc}
\immediate\write\gtoutfile{MSC-class: \theprimaryclass\ifx\thesecondaryclass\relax\else, \thesecondaryclass\fi}
\immediate\write\gtoutfile{Journal-ref: Algebr. Geom. Topol. \thevolumenumber\s
(\thevolumeyear) \startpage-\finishpage}
\immediate\write\gtoutfile{Comments: Published by Algebraic and
Geometric Topology at}
\immediate\write\gtoutfile{\s\s\s  http://www.maths.warwick.ac.uk/agt/AGTVol\thevolumenumber/agt-\thevolumenumber-\thepapernumber.abs.html}
\immediate\write\gtoutfile{\noexpand\\}
\immediate\write\gtoutfile{}
\ifx\theasciiabstract\relax
\immediate\write\gtoutfile{\theabstract}\else
\immediate\write\gtoutfile{\theasciiabstract}\fi
\immediate\write\gtoutfile{}
\immediate\write\gtoutfile{\noexpand\\}
\immediate\write\gtoutfile{}
\immediate\closeout\gtoutfile}}  

\def\maketitlepage{\makeagttitle\makeheadfile}
\let\makeshorttitle\maketitlepage
\let\maketitle\maketitlepage

\lognumber{16}
\volumenumber{5}
\volumeyear{2005}
\papernumber{16}
\pagenumbers{369}{378}
\received{10 March 2005} 
\revised{25 March 2005}
\accepted{12 April 2005}
\published{30 April 2005}

\usepackage{amsmath,amssymb,graphicx}

\newtheorem{LM}{Lemma}[section]
\newtheorem{THM}[LM]{Theorem}

\newtheorem{CO}[LM]{Corollary}
\newtheorem{CL}[LM]{Claim}

\theoremstyle{remark}
\newtheorem*{rem}{Remark}
\newtheorem*{rems}{Remarks}

\begin{document}

\title
{All integral slopes can be Seifert fibered slopes\\for hyperbolic knots} 

\author{Kimihiko Motegi\\Hyun-Jong Song}
\address{Department of Mathematics, Nihon University\\Tokyo 156-8550, Japan}
\secondaddress{Division of Mathematical Sciences, Pukyong National 
University\\599-1 Daeyondong, Namgu, Pusan 608-737, Korea}
\asciiaddress{Department of Mathematics, Nihon University\\Tokyo 
156-8550, Japan\\and\\Division of Mathematical Sciences, 
Pukyong National University\\599-1 Daeyondong, Namgu, Pusan 608-737, Korea}

\gtemail{\mailto{motegi@math.chs.nihon-u.ac.jp}{\rm\qua
and\qua}\mailto{hjsong@pknu.ac.kr}}
\asciiemail{motegi@math.chs.nihon-u.ac.jp, hjsong@pknu.ac.kr}

\begin{abstract}
	Which slopes can or cannot appear as Seifert fibered slopes 
for hyperbolic knots in the $3$-sphere $S^3$? 
It is conjectured that if $r$-surgery on a hyperbolic knot in 
$S^3$ yields a Seifert fiber space, 
then $r$ is an integer. 
	We show that for each integer $n \in \mathbb{Z}$, 
there exists a tunnel number one, 
hyperbolic knot $K_n$ in $S^3$ such that 
$n$-surgery on $K_n$ produces a small Seifert fiber space. 
\end{abstract}
\asciiabstract{%
Which slopes can or cannot appear as Seifert fibered slopes for
hyperbolic knots in the 3-sphere S^3?  It is conjectured that if
r-surgery on a hyperbolic knot in S^3 yields a Seifert fiber space,
then r is an integer.  We show that for each integer n, there exists a
tunnel number one, hyperbolic knot K_n in S^3 such that n-surgery on
K_n produces a small Seifert fiber space.}


\maketitle
{\small\it This paper is dedicated to Donald M. Davis on the occasion of his
60th birthday.\leftskip25pt\rightskip25pt\par}

\section{Introduction}
	Let $K$ be a knot in the $3$-sphere $S^3$ 
with a tubular neighborhood $N(K)$. 
Then the set of \textit{slopes} for $K$ 
(i.e., $\partial N(K)$-isotopy classes of simple loops on $\partial N(K)$) 
is identified with $\mathbb{Q} \cup \{\infty \}$ 
using preferred meridian-longitude pair so that 
a meridian corresponds to $\infty$. 
A slope $\gamma$ is said to be \textit{integral} if 
a representative of $\gamma$ intersects a meridian exactly once, 
in other words, $\gamma$ corresponds to an integer under the above 
identification. 
In the following, 
we denote by $(K; \gamma)$ 
the $3$-manifold obtained from $S^3$ by Dehn surgery
on a knot $K$ with slope $\gamma$, 
i.e., by attaching a solid torus to $S^3-$int$N(K)$ in such a way 
that $\gamma$ bounds a meridian disk of the filled solid torus. 
If $\gamma$ corresponds to $r \in \mathbb{Q} \cup \{ \infty \}$, 
then we identify $\gamma$ and $r$ and write $(K; r)$ for $(K; \gamma)$. \par

	We denote by $\mathcal{L}$ 
the \textit{set of lens slopes} 
$\{r \in \mathbb{Q}\ |\ \exists$ 
hyperbolic knot $K \subset S^3$ such that $(K;r)$ is a lens space$\}$, 
where $S^3$ and $S^2 \times S^1$ are also considered as lens spaces. 
Then the cyclic surgery theorem \cite{CGLS} implies that 
$\mathcal{L} \subset \mathbb{Z}$. 
A result of Gabai \cite[Corollary 8.3]{Ga} shows that 
$0 \not\in \mathcal{L}$, 
a result of Gordon and Luecke \cite{G-Lu} shows 
that $\pm 1 \not\in \mathcal{L}$. 
In \cite{KM} Kronheimer and Mrowka 
prove that $\pm 2 \not\in \mathcal{L}$. 
Furthermore, 
a result of Kronheimer, Mrowka, Ozsv\'ath and Szab\'o \cite{KMOS} 
implies that $\pm 3, \pm 4 \not\in \mathcal{L}$.  
Besides, 
Berge \cite[Table of Lens Spaces]{Berge} 
suggests that 
if $n \in \mathcal{L}$, 
then $|n| \ge 18$ and 
not every integer $n$ with $|n| \ge 18$ appears  
in $\mathcal{L}$. 
Fintushel and Stern \cite{FS} had shown that $18$-surgery on 
the $(-2, 3, 7)$ pretzel knot yields a lens space. 

\textit{Which slope $($rational number$)$ can or cannot appear in 
the set of Seifert fibered slopes 
$\mathcal{S} = \{r \in \mathbb{Q}\ |\ \exists$ 
hyperbolic knot $K \subset S^3$ such that $(K;r)$ 
is Seifert fibered$\}$?} 
It is conjectured that $\mathcal{S} \subset \mathbb{Z}$ \cite{Go}. 

The purpose of this paper is to prove: 

\begin{THM}
\label{Seifert slopes}
For each integer $n \in \mathbb{Z}$, 
there exists a tunnel number one, 
hyperbolic knot $K_n$ in $S^3$ 
such that $(K_n; n)$ is a small Seifert fiber space 
$($i.e., a Seifert fiber space over $S^2$ 
with exactly three exceptional fibers$)$. 
\end{THM}

\begin{rem}
Since $K_n$ has tunnel number one, 
it is embedded in a genus two Heegaard surface of $S^3$ and 
strongly invertible \cite[Lemma 5]{Mor}. 
See \cite[Question 3.1]{MMM}. 
\end{rem}

Theorem \ref{Seifert slopes}, 
together with the previous known results, 
shows: 

\begin{CO}
\label{LZS}
$\mathcal{L} \subsetneqq \mathbb{Z} \subset \mathcal{S}$. 
\end{CO}

\begin{rems}$\phantom{99}$

(1)\qua For the \textit{set of reducing slopes} 
$\mathcal{R} = \{r \in \mathbb{Q}\ |\ \exists$ 
hyperbolic knot $K \subset S^3$ such that $(K;r)$ is reducible$\}$, 
Gordon and Luecke \cite{G-Lu1} have shown that 
$\mathcal{R} \subset \mathbb{Z}$. 
In fact, the cabling conjecture \cite{GS} asserts that 
$\mathcal{R} = \emptyset$. 

(2)\qua For the \textit{set of toroidal slopes} 
$\mathcal{T} = \{r \in \mathbb{Q}\ |\ \exists$ 
hyperbolic knot $K \subset S^3$ such that $(K;r)$ is toroidal$\}$, 
Gordon and Luecke \cite{G-Lu2} have shown that 
$\mathcal{T} \subset \mathbb{Z}/2$ (integers or half integers). 
In \cite{Tera}, 
Teragaito shows that 
$\mathbb{Z} \subset \mathcal{T}$ and conjectures that 
$\mathcal{T} \subsetneqq \mathbb{Z}/2$. 
\end{rems}

\textbf{Acknowledgements}\qua
We would like to thank the referee for careful reading 
and useful comments. \newline
The first author was partially 
supported by Grant-in-Aid for 
Scientific Research (No.\ 15540095), 
The Ministry of Education, Culture, Sports, 
Science and Technology, Japan. \newline

\section{Hyperbolic knots with Seifert fibered surgeries}

	Our construction is based on an example of 
a longitudinal Seifert fibered surgery given in \cite{IMS}. 

\begin{figure}[ht!]
\begin{center}
\includegraphics[width=0.4\linewidth]{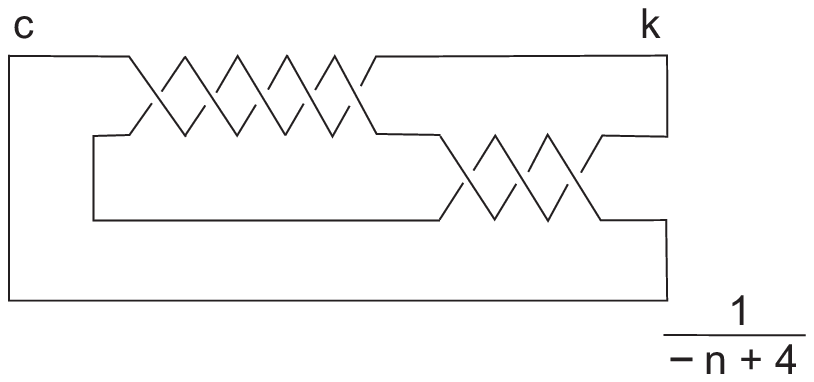}
\caption{$K_n$}
\label{fig:knotKn}
\end{center}
\end{figure}

Let $k \cup c$ be a $2$-bridge link given in 
Figure \ref{fig:knotKn}, 
and let $K_n$ be a knot obtained from $k$ by 
$\frac{1}{-n+4}$-surgery along $c$. 

        We shall say that a Seifert fiber space is of 
\textit{type} $S^2(n_1, n_2, n_3)$ 
if it has a Seifert fibration over 
$S^2$ with three exceptional fibers of 
indices $n_1, n_2$ and $n_3$ $(n_i \ge 2)$.  
Since $K_4$ is unknotted, 
$(K_4; 4)$ is a lens space $L(4, 1)$. 
For the other $n$'s, 
we have: 

\begin{LM}
\label{Seifert}
$(K_n; n)$ is a small Seifert fiber space of type 
$S^2(3, 5, |4n-15|)$ for any integer $n \ne 4$. 
\end{LM}

\begin{proof}
	Since the linking number of 
$k$ and $c$ is one (with suitable orientations), 
$(K_n ; n)$ has surgery descriptions 
as in Figure \ref{fig:description}. 

\begin{figure}[ht!]
\begin{center}
\includegraphics[width=1.0\linewidth]{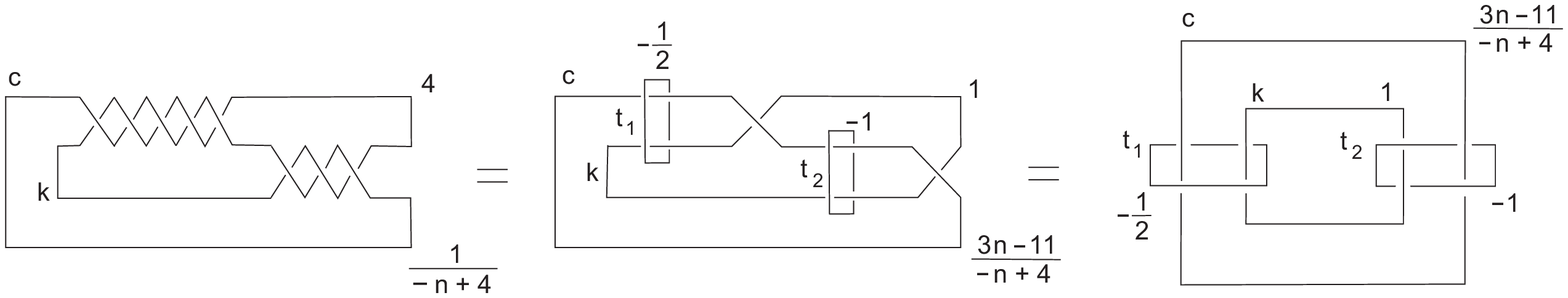}
\caption{Surgery descriptions of $(K_n; n)$}
\label{fig:description}
\end{center}
\end{figure}

	Let us take the quotient by the strong inversion of $S^3$ 
with an axis $L$ as shown in Figure \ref{fig:msequence1}.

\begin{figure}[ht!]
\begin{center}
\includegraphics[width=0.9\linewidth]{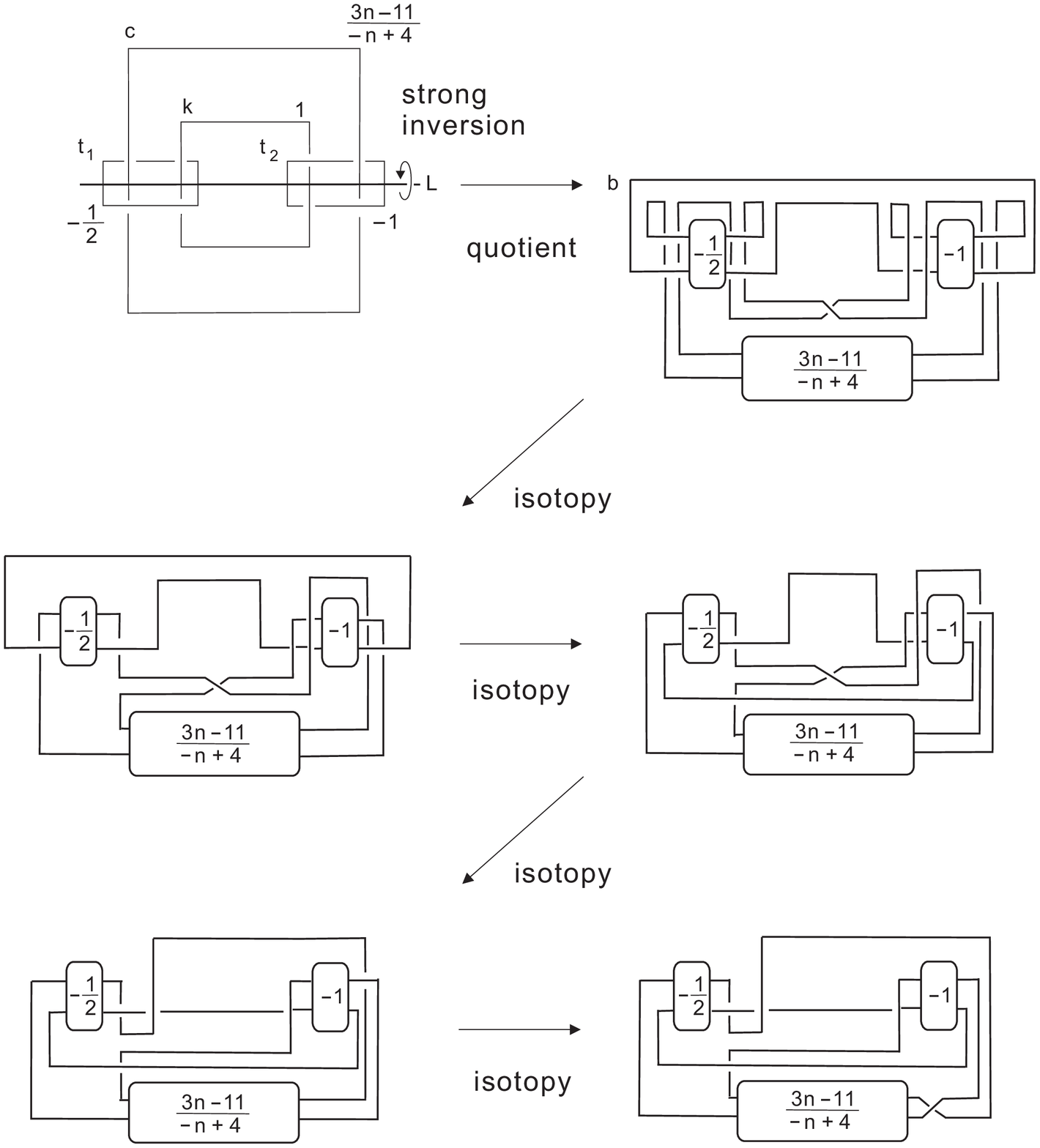}
\nocolon\caption{}
\label{fig:msequence1}
\end{center}
\end{figure}

Then we obtain a branch knot $b'$ which is the image of the axis $L$. 
The Montesinos trick (\cite{Mon}, \cite{Bl}) shows that 
$-\frac{1}{2}, -1, \frac{3n-11}{-n+4}$ and 
$1$-surgery on 
$t_1, t_2, c$ and $k$ in the upstairs correspond to 
$-\frac{1}{2}, -1, \frac{3n-11}{-n+4}$ and 
$1$-untangle surgery on $b'$ in the downstairs, 
where an \textit{$r$-untangle surgery} is 
a replacement of $\frac{1}{0}$-untangle by 
$r$-untangle. 
(We adopt Bleiler's convention \cite{Ble} 
on the parametrization of rational tangles.) 
These untangle surgeries convert $b'$ into 
a link $b$ (Figure \ref{fig:msequence1}).

\begin{figure}[ht!]
\begin{center}
\includegraphics[width=0.9\linewidth]{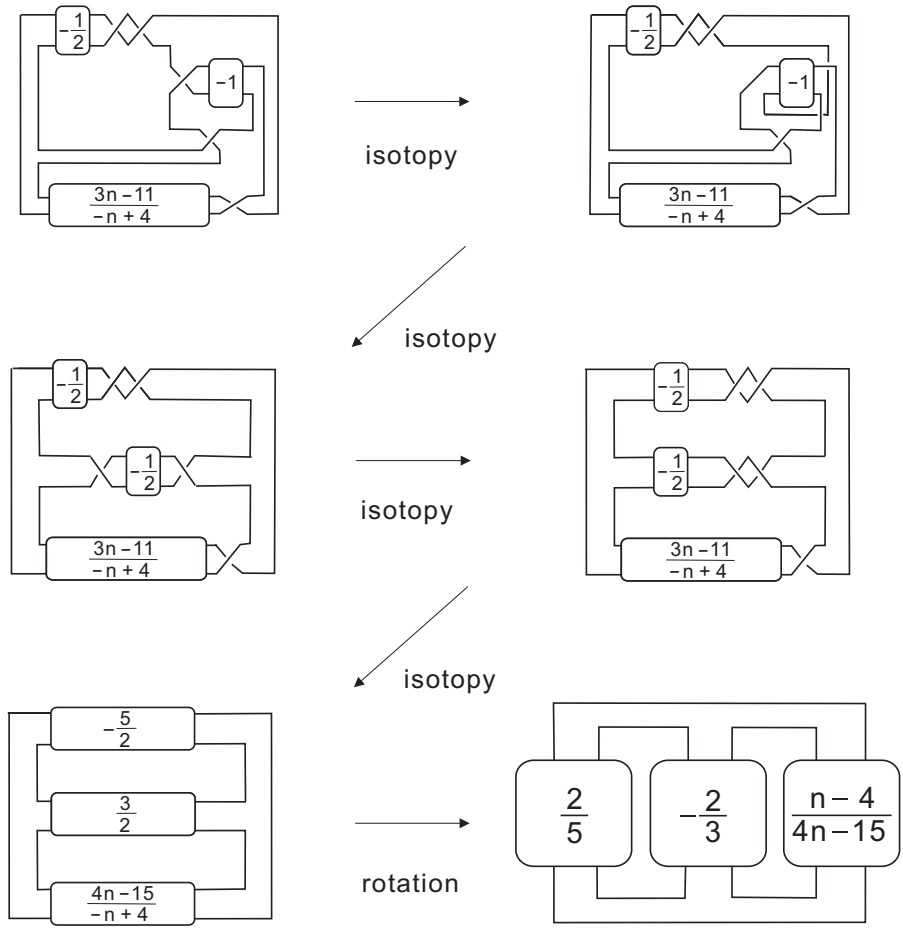}
\caption{Continued from Figure \ref{fig:msequence1}}
\label{fig:msequence2}
\end{center}
\end{figure}

	Following the sequence of isotopies in 
Figures \ref{fig:msequence1} and \ref{fig:msequence2}, 
we obtain a Montesinos link 
$M(\frac{2}{5}, -\frac{2}{3}, \frac{n-4}{4n-15})$. 

Since $(K_n ; n)$ is the double branched cover of $S^3$ branched over 
the Montesinos link $M(\frac{2}{5}, -\frac{2}{3}, \frac{n-4}{4n-15})$, 
$(K_n; n)$ is a Seifert fiber space of type 
$S^2(3, 5, |4n-15|)$ as desired. 
\end{proof}

\begin{LM}
\label{hyperbolic}
The knot $K_n$ is hyperbolic if 
$n \ne 3, 4, 5$. 
\end{LM}

\begin{proof}
Note that the $2$-bridge link given in Figure \ref{fig:knotKn} 
is not a $(2, p)$-torus link, 
and hence by \cite{Me} 
it is a hyperbolic link. 
If $n \ne 3, 4, 5$, 
then $|-n+4|>1$ and 
it follows from 
\cite[Theorem 1]{AMM1} (also \cite[Theorem 1.2]{AMM3})
that $K_n$ is a hyperbolic knot. 
See also 
\cite[Corollary A.2]{G-Lu3}, 
\cite[Theorem 1.2]{MM3} and 
\cite[Theorem 1.1]{AMM2}. 
\end{proof}

\begin{rem}
It follows from \cite{Ma}, \cite{KMS} 
that $K_n$ is a nontrivial knot except when $n = 4$. 
An experiment using Weeks' computer program 
``SnapPea" \cite{W} suggests that 
$K_3$ and $K_5$ are hyperbolic, 
but we will not use this experimental results. 
\end{rem}

\begin{LM}
\label{tunnel}
The knot $K_n$ has tunnel number one for any integer $n \ne 4$. 
\end{LM}

\begin{proof}
Since the link $k \cup c$ is a two-bridge link, 
the tunnel number of $k \cup c$ is one 
with unknotting tunnel $\tau$; 
A regular neighborhood $N(k \cup c \cup \tau)$ is a genus two 
handlebody and $S^3 - \mathrm{int}N(k \cup c \cup \tau)$ is also 
a genus two handlebody, 
see Figure \ref{fig:tunnel}. 

\begin{figure}[ht!]
\begin{center}
\includegraphics[width=0.5\linewidth]{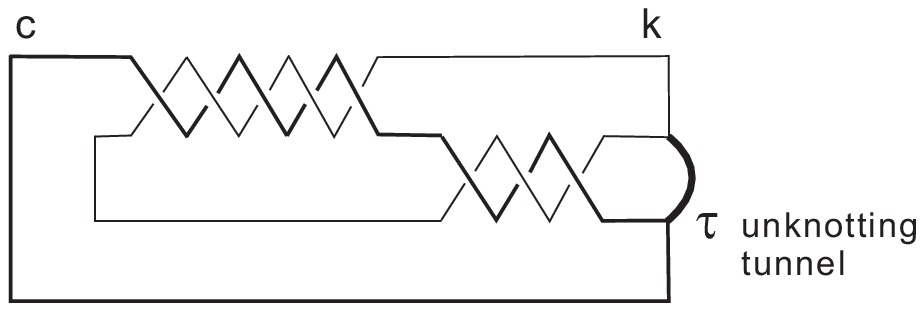}
\nocolon\caption{}
\label{fig:tunnel}
\end{center}
\end{figure}

Then the general fact below 
(in which $k \cup c$ is not necessarily a two-bridge link) 
shows that the tunnel number of $K_n$ is less than or equal to one. 
Since our knot $K_n$ $(n \ne 4)$ is knotted in $S^3$, 
the tunnel number of $K_n$ is one.   
\end{proof}

\begin{CL}
\label{tunnel1}
Let $k \cup c$ be a two component link in $S^3$ 
which has tunnel number one. 
Assume that $c$ is unknotted in $S^3$. 
Then every knot obtained from $k$ by twisting along $c$ 
has tunnel number at most one. 
\end{CL}

\begin{proof}
	Let $\tau$ be an unknotting tunnel and 
$V$ a regular neighborhood of $k \cup c \cup \tau$ in $S^3$; 
$V$ is a genus two handlebody. 
Since $\tau$ is an unknotting tunnel for $k \cup c$, 
by definition, 
$W = S^3 - \mathrm{int}V$ is also 
a genus two handlebody. 
Take a small tubular neighborhood $N(c) \subset \mathrm{int}V$ 
and perform $-\frac{1}{n}$-surgery on $c$ using $N(c)$. 
Then we obtain 
a knot $k_n$ as the image of $k$ and 
obtain a genus two handlebody $V(c; -\frac{1}{n})$. 
Note that $V(c; -\frac{1}{n})$ and $W$ define 
a genus two Heegaard splitting of $S^3$, 
see Figure \ref{fig:arc}, 
where $c_n^*$ denotes the core of the filled solid torus. 

\begin{figure}[ht!]
\begin{center}
\includegraphics[width=0.88\linewidth]{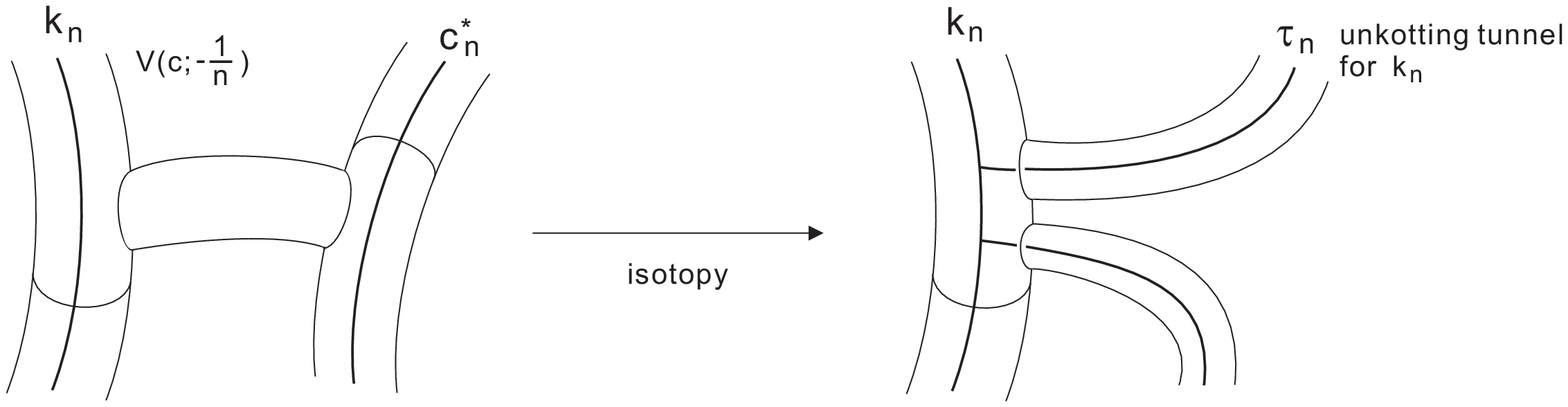}
\nocolon\caption{}
\label{fig:arc}
\end{center}
\end{figure}

Then it is easy to see that an arc $\tau_n$ given by Figure \ref{fig:arc} 
is an unknotting tunnel for $k_n$ as desired. 
\end{proof}

	Now we are ready to prove Theorem \ref{Seifert slopes}. 
Lemmas \ref{Seifert}, \ref{hyperbolic} and \ref{tunnel} 
show that our knots $K_n$ enjoy the required properties, 
except for $n = 3, 4, 5$. 
To prove Theorem \ref{Seifert slopes}, 
we find hyperbolic knots $K'_n$ so that 
$(K'_n; n)$ is Seifert fibered for $n = 3, 4, 5$ 
(instead of showing that $K_3$, $K_5$ are hyperbolic). 
As the simplest way, 
let $K'_3$, $K'_4$ and $K'_5$ 
be the mirror image of 
$K_{-3}$, $K_{-4}$ and $K_{-5}$, 
respectively. 
Since $K_{-3}$, $K_{-4}$ and $K_{-5}$ are tunnel number one, 
hyperbolic knots 
by Lemmas \ref{hyperbolic} and \ref{tunnel}, 
their mirror images 
$K'_3$, $K'_4$ and $K'_5$ are also tunnel number one, 
hyperbolic knots. 
It is easy to observe that 
$(K'_3; 3)$ 
(resp. $(K'_4; 4)$, $(K'_5; 5)$) 
is the mirror image of $(K_{-3}; -3)$ 
(resp. $(K_{-4}; -4)$, $(K_{-5}; -5)$). 
By Lemma \ref{Seifert}, 
$(K_{-3}; -3)$, $(K_{-4}; -4)$ and $(K_{-5}; -5)$ are Seifert fibered, 
and hence $(K'_3; 3)$, $(K'_4; 4)$ and $(K'_5; 5)$ are also Seifert fibered. 
Putting $K_n$ as $K'_n$ for $n = 3, 4, 5$, 
we finish a proof of Theorem \ref{Seifert slopes}. 
\hspace*{\fill} $\qed$

\section{Identifying exceptional fibers}
	In \cite{MM3}, 
Miyazaki and Motegi conjectured that 
if $K$ admits a Seifert fibered surgery, 
then there is a trivial knot $c \subset S^3$ disjoint from $K$ 
which becomes a Seifert fiber in the resulting Seifert fiber space, 
and verified the conjecture for several Seifert fibered surgeries 
\cite[Section 6]{MM3}, see also \cite{EM}. 
Furthermore, 
computer experiments via ``SnapPea" \cite{W} suggest that 
such a knot $c$ is realized by a short closed geodesic in 
the hyperbolic manifold $S^3 - K$,
for details see \cite[Section 9]{MM3}, \cite{Mot2}. 

	In this section, 
we verify the conjecture for Seifert fibered surgeries 
given in Theorem \ref{Seifert slopes}. 

	Recall that $K_n$ is obtained from $k$ by 
$\frac{1}{-n+4}$-surgery on the trivial knot $c$ 
(i.e., $(n-4)$-twist along $c$), see Figure \ref{fig:knotKn}. 
Denote by $c_n$ the core of the filled solid torus. 
Then $K_n \cup c_n$ is a link in $S^3$ 
such that 
$c_n$ is a trivial knot. 

\begin{LM}
\label{fiber}
After $n$-surgery on $K_n$, 
$c_n$ becomes an exceptional fiber of index $|4n-15|$ 
in the resulting Seifert fiber space $(K_n; n)$. 
\end{LM}

\begin{proof}
	Following the sequences given by 
Figures \ref{fig:msequence1} and \ref{fig:msequence2}, 
we have a Montesinos link with three arcs
$\gamma$, $\tau_1$ and $\tau_2$ as in Figure \ref{fig:positions}, 
where $n = 1$ in the final Montesinos link, 
and $\gamma$, $\tau_1$, $\tau_2$ and $\kappa$ are the 
images of $c$, $t_1$, $t_2$ and $k$, respectively. 

\begin{figure}[ht!]
\begin{center}
\includegraphics[width=1.0\linewidth]{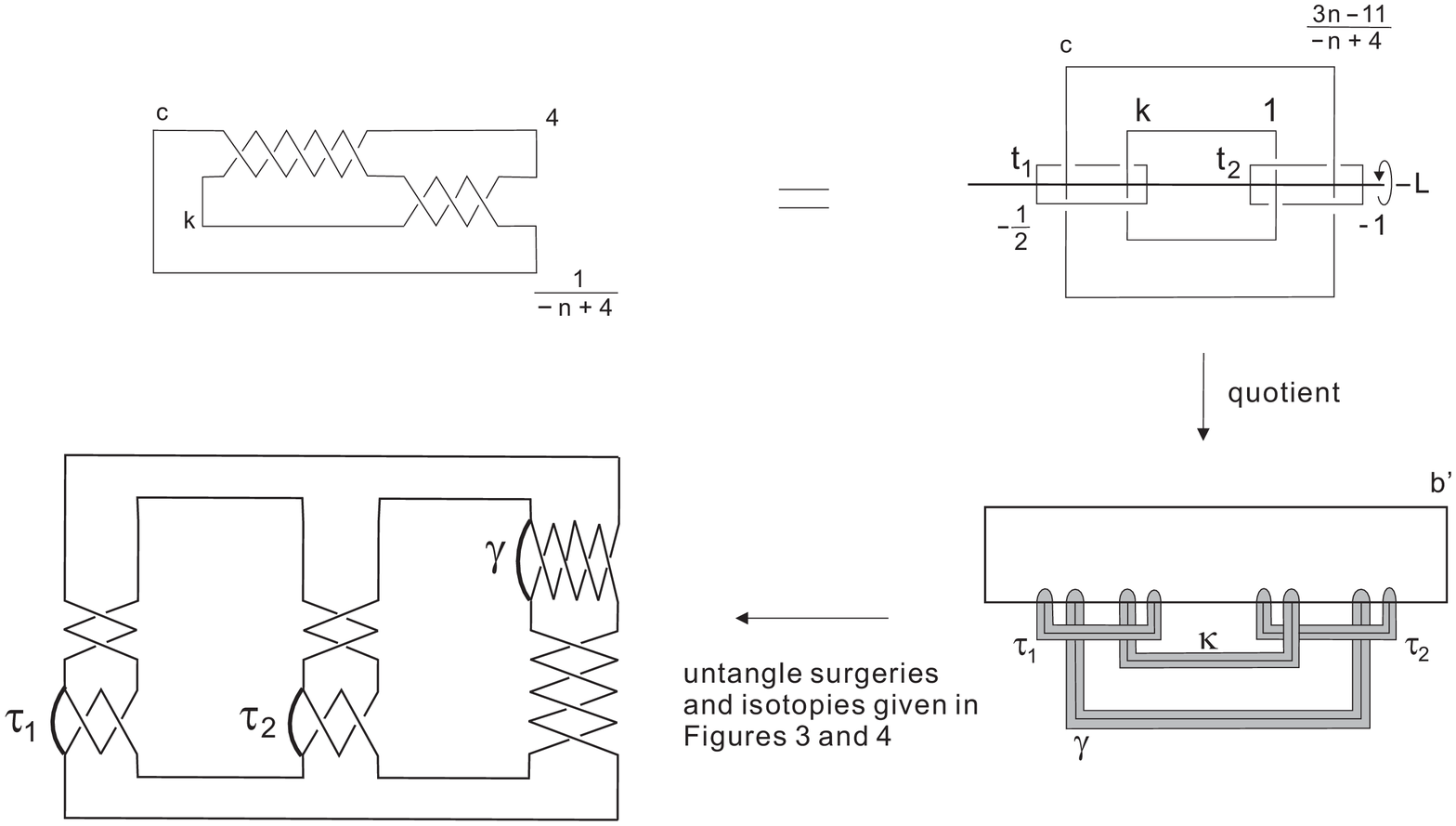}
\caption{Positions of exceptional fibers}
\label{fig:positions}
\end{center}
\end{figure}

From Figure \ref{fig:positions} we recognize that 
$t_1, t_2$ and $c$ become exceptional fibers of 
indices $5$, $3$ and $|4n-15|$, respectively 
in $(K_n; n)$. 
\end{proof}

	For $n \ne 3, 4, 5$, 
$c_n$ becomes an exceptional fiber of index $|4n-15|$, 
which is the unique maximal index,  
in $(K_n; n)$. 
	Experiments via ``SnapPea" \cite{W} 
suggest that $c_n$ is a shortest closed geodesic in $S^3 - K_n$ $(n \ne 3, 4, 5)$. 
For sufficiently large $|n|$, 
hyperbolic Dehn surgery theorem \cite{T1}, \cite{T2} 
shows that 
$c_n$ is the unique shortest closed geodesic in $S^3 - K_n$. \par

	Let us assume that $n = 3, 4, 5$. 
Then we have put $K_n$ as the mirror image of $K_{-n}$ in 
the proof of Theorem \ref{Seifert slopes}. 
Let $k' \cup c'$ be the mirror image of the link $k \cup c$. 
Then $K_n$ is obtained also from $k'$ by 
$\frac{1}{-n-4}$-surgery on $c'$ (i.e., $(n+4)$-twist along $c'$); 
we denote the core of the filled solid torus by $c'_n$.  
Note that there is an orientation reversing diffeomorphism 
from $(K_{-n}; -n)$ to $(K_n; n)$ sending 
$c_{-n}$ (regarded as a fiber in $(K_{-n}; -n)$) 
to $c'_n$ (regarded as a fiber in $(K_n; n)$). 
Thus the above observation implies that 
$c'_n$ becomes 
an exceptional fiber of index $|4n+15|$, 
which is the unique maximal index, 
in $(K_n; n)$ $(n = 3, 4, 5)$.

\Addresses\recd

\end{document}